%% AMS-TEX

\input amstex
\documentstyle{amsppt}

\define\a{\alpha} 

\predefine\barund{\b}
\redefine\b{\beta}
\define\brel{\buildrel}
\define\bsk{\bigskip}
\define\bu{\bullet}
\define\CC{{\Bbb C}}
\define\cA{{\Cal A}}

\define\cE{{\Cal E}}
\define\cF{{\Cal F}}
\define\cG{{\Cal G}}

\define\cl{\colon}
\define\cL{{\Cal L}}
\define\cM{{\Cal M}}

\define\cO{{\Cal O}}
\define\cod{\operatorname{cod}}

\define\cQ{{\Cal Q}}

\define\cU{{\Cal U}}

\define\cZ{{\Cal Z}}
\predefine\dotund{\d}
\redefine\d{\delta}
\define\del{\partial}
\predefine\dotov{\D}
\redefine\D{\Delta}
\define\Def{\operatorname{Def}}

\define\e{\epsilon}

\define\End{\operatorname{End}}
\define\Ext{\operatorname{Ext}}
\define\es{\emptyset}

\define\G{\Gamma}

\define\Hom{\operatorname{Hom}}
\define\hra{\hookrightarrow}

\predefine\idot{\i}
\redefine\i{\iota}

\predefine\Imaginary{\Im}
\redefine\Im{\operatorname{Im}}
\define\k{\kappa}
\define\Ker{\operatorname{Ker}}
\predefine\lcut{\l}
\redefine\l{\lambda}

\predefine\Lcut{\L}
\redefine\L{\Lambda}
\define\lra{\longrightarrow}
\define\msk{\medskip}
\define\n{\noindent}
\define\NN{{\Bbb N}}

\predefine\ocut{\o}
\redefine\o{\omega}
\define\op{\oplus}
\define\ot{\otimes}
\define\ov{\overline}
\predefine\Ocut{\O}
\redefine\O{\Omega}
\define\PP{{\Bbb P}}
\define\Pic{\operatorname{Pic}}
\define\QQ{{\Bbb Q}}

\define\rk{\operatorname{rk}}

\define\s{\sigma}
\predefine\Spar{\S}
\redefine\S{\Sigma}

\predefine\sgerman{\ss}
\redefine\ss{\subset}

\predefine\tie{\t}
\redefine\t{\theta}
\define\T{\Theta}

\define\tm{\times}

\define\Tr{\operatorname{Tr}}

\define\ul{\underline}
\define\vf{\varphi}

\define\wh{\widehat}
\define\wt{\widetilde}

\define\ZZ{{\Bbb Z}}

\define\mx{\cM_{\xi}}
\define\dx{\D_{\xi}}

\topmatter
\title
Moduli of Vector-Bundles on Surfaces. 
\endtitle
\rightheadtext{Moduli of Vector-Bundles on Surfaces}
\author 
Kieran G.~O'Grady
\endauthor
\address
Universit\`a di Salerno, Facolt\`a di Scienze, Baronissi (Sa) - Italia
\endaddress
\email
ogrady\@ mat.uniroma1.it
\endemail
\thanks 
Partially supported by G.N.S.A.G.A. (C.N.R.).
\endthanks
\date 
September 14 1996
\enddate
\endtopmatter
\document
\head
0. Introduction.
\endhead
In the 1980's Donaldson proved some  spectacular new results on 
classification of $C^{\infty}$ four-manifolds by studying 
anti-self-dual (ASD) connections on an  $SU(2)$-bundle.  If the
four-manifold underlies a complex projective surface, the set of ASD
connections modulo gauge transformations is identified with a
moduli space of slope-stable vector-bundles on the surface.
Donaldson~\cite{D, \Spar V} proved that for rank two these moduli spaces
are generically smooth of the expected dimension (see Section~(1) for
precise definitions), provided the expected dimension is large enough;  this
implies that the polynomial invariants of a projective surface are not 
zero. In this paper we will present algebro-geometric results which were
inspired  by  Donaldson's theory;  there will be  no discussion of
relations with Gauge theory. First of all we will sketch our proof~\cite{O2}
of a  theorem  proved also by Gieseker-Li~\cite{GL1,GL2}.
\proclaim{(0.1) Theorem (Gieseker-Li, O'Grady)}
Let $S$ be an irreducible  smooth complex projective surface, and $H$ an
ample divisor on $S$. There exists $\D(r)$ such that the moduli of
$H$-semistable (in the sense of Gieseker-Maruyama) rank-$r$ torsion-free
sheaves on $S$, with Chern classes $c_1,c_2\in H^*(S;\ZZ)$, is reduced of
the expected dimension 
$$2rc_2-(r-1)c_1^2-(r^2-1)\chi(\cO_S)+h^1(\cO_S),$$
provided  $\D:=c_2-\left((r-1)/2r\right)c_1^2>\D(r)$. Furthermore (for
$\D\gg 0$) the open subset parametrizing $H$-slope-stable vector-bundles
is   dense, and the moduli space is irreducible. 
\endproclaim
The above  statement requires  a few comments. If $r=1$ the moduli space
$\cM$  is isomorphic to the product of $Pic^{c_1}(S)$ and the Hilbert scheme
parametrizing length-$c_2$ zero-dimensional subschemes of $S$: as is
well-known this Hilbert scheme is always smooth, irreducible and of the
expected dimension, hence so is $\cM$. We are really concerned with the
case $r\ge 2$: from now on we will always assume that the rank is at least
two, unless we specify otherwise. We deal with Gieseker-Maruyama 
semistable torsion-free sheaves, rather than with slope-stable
vector-bundles, because of a theorem of Gieseker and
Maruyama~\cite{G1,Ma}: The moduli space of  semistable torsion-free
sheaves (containing the moduli space of slope-stable vector-bundles as an
open subscheme) is projective. Regarding the hypothesis that $\D$ is large:
The moduli space is empty if $\D<0$, by Bogomolov's Inequality, and on the
other hand it is non-empty  if $\D\gg 0$~\cite{Ma,LQ,G2}. For "low"
non-negative values of $\D$ there are many examples~\cite{G2,O2} of
moduli spaces which are not of the expected dimension (or which are
reducible~\cite{Me}): this is a tipical phenomenon occurring for surfaces of
Kodaira dimension at least one. At the other extreme of the
Kodaira-Enriques classification, if say $S$ is the projective plane, then
Theorem~(0.1) has been known for a long time in a stronger
form~\cite{Ma,DL}. More  generally, if $S$ is not of general type the  moduli
space can be somewhat analyzed~\cite{Ba,H,ES,Mk1,Mk2,F1} because its
structure reflects the special properties of $S$ given by the
Enriques-Kodaira classification. If instead $S$ is of general type, very little
is known about  moduli of vector-bundles; Theorem~(0.1) is one of the few
general results. After sketching a proof of this theorem we will discuss
holomorphic two-forms on the moduli space (of sheaves with fixed
determinant). There is a natural map,  first studied by Mukai, associating to
a holomorphic two-form $\o$ on $S$ a holomorphic two-form $\o_{\xi}$ on
the moduli space. If the rank is two  and some other hypotheses are
satisfied, then $\o_{\xi}$ is non-degenerate at the generic
point~\cite{Mk1,O1}. As noticed by Tyurin~\cite{T} the non-degeneracy of
$\o_{\xi}$ implies that the image of the map       
$$\matrix
\text{moduli space} & \to & CH_0(S) \\
[E] & \mapsto & c_2(E)
\endmatrix$$
has "dimension" equal to that of the moduli space. Finally we will discuss
the  Kodaira dimension of the moduli space. We will sketch J.~Li's
proof~\cite{L2} that if $S$ is of general type then the moduli space  is of
general type, if the rank is two and  certain other hypotheses are satisfied.
The results on two-forms and the Kodaira dimension  had been proved when
Theorem~(0.1) was known in rank two only. We observe that
 since~(0.1) holds in arbitrary rank,  analogous
results on the non-degeneracy of $\o_{\xi}$, and on the Kodaira dimension
of the moduli space, are valid if a certain conjecture~(2.4) regarding
vector-bundles on curves is true. We will verify this conjecture for
arbitrary rank and a special choice of degree (Proposition~(2.5)). 
\subhead
Notation
\endsubhead 
All schemes are defined over $\CC$. We let $S$ be a smooth irreducible
projective surface, and  $K$  be its canonical divisor class. We let $H$
be an ample divisor  on $S$. 
\msk
\n
Let $X$ be a projective variety of dimension $n$, and $D$ be an ample
divisor on $X$: for a torsion-free sheaf $F$ on $X$ one sets
$$\text{slope of $F$}=\mu(F):={c_1(F)\cdot D^{n-1}\over \rk(F)},
\qquad p_F(n):={\chi\left(F\ot\cO_X(nD)\right)\over \rk(F)}.$$
The sheaf $F$ is $D$-{\it slope-semistable} (respectively $D$-{\it
semistable}) if
$$\mu(E)\le\mu(F)\qquad \text{($p_E(n)\le p_F(n)$ for  $n\gg 0$),}$$
for all (non-zero) subsheaves $E\ss F$; if strict inequality holds whenever
$\rk(E)<\rk(F)$ then $F$ is $D$-{\it slope-stable} (respectively   $D$-{\it
stable}). One easily checks the implications:
$$D-\text{slope-stable}\Longrightarrow D-\text{stable} \Longrightarrow
D-\text{semistable}\Longrightarrow D-\text{slope-semistable}.$$
Now let's specialize to the case $X=S$. For a torsion-free sheaf
$F$ on $S$ the discriminant is  
$$\D_F:=c_2(F)-{\rk(F)-1\over 2\rk(F)}c_1(F)^2.$$
We label moduli spaces of sheaves on $S$ with triples of {\it sheaf data}
$$\xi=(\rk(\xi),\det(\xi),c_2(\xi))\in \NN\tm \Pic(S)\tm H^4(S;\ZZ),$$
and we set
$$\displaylines{
\quad\mx(S,H):=\{\text{$H$-s.s.~tors.-free sheaf $F$ on $S$}\hfill\cr
\hfill \text{with $\rk(F)=\rk(\xi)$, $\det F\cong\det(\xi)$,
$c_2(F)=c_2(\xi)$}\}/\text{S-equivalence}.\quad\cr}$$ 
To define S-equivalence one considers a Jordan-H\"older(JH) filtration 
$$0=F_0\ss F_1\ss \cdots \ss F_n=F,$$
i.e.~such that $p_{F_i}=p_F$ and $F_i/F_{i-1}$ is stable for
$i=1,\ldots,n$. The associated graded sheaf $Gr_{JH}(F):=\bigoplus_{i=1}^n 
F_i/F_{i-1}$ is unique up to isomorphism (although a JH filtration is not
unique): two semistable sheaves $F$, $F'$ are {\it $S$-equivalent} if
$Gr_{JH}(F)\cong Gr_{JH}(F')$. Thus $\mx(S,H)$ contains an open subscheme
$\mx^{st}(S,H)$ parametrizing $\ul{\text{isomorphism}}$ classes of stable
sheaves. By a theorem of Gieseker and Maruyama~\cite{G1,Ma}, $\mx(S,H)$ is
projective.  We indicate by $[F]$ the point of $\mx(S,H)$ corresponding to a
semistable sheaf $F$.  We set  
$$c_1(\xi):=c_1\left(\det(\xi)\right)\qquad
\D_{\xi}:=c_2(\xi)-{\rk(\xi)-1\over 2\rk(\xi)}c_1(\xi)^2.$$
Notice that we fix
the determinant of  sheaves, not just $c_1\in H^2$ as in Theorem~(0.1).
\remark{\rom{(0.2)} Remark}
How does the moduli space vary when we change
the polarization $H$? This problem is studied in various papers (for
example~\cite{Q,MW}). We will not discuss the known results, except for the
following general fact. Let $H_1$, $H_2$ be ample divisors on $S$, and fix
the rank of the sheaves: if $\dx$ is sufficiently large the moduli spaces
$\mx(S,H_1)$, $\mx(S,H_2)$ are birational. Thus for many purposes we can
fix the polarization $H$, and this is what we will always do. To simplify
notation we write $\mx$ instead of $\mx(S,H)$.
\endremark   
\n
A {\it family of sheaves on $X$ parametrized by $B$} consists of a sheaf on
$X\tm B$, flat over $B$. We say $\mx$ is a {\it fine moduli space} if
$\mx^{st}=\mx$ (i.e.~semistablity implies
stability),  and furthermore there exists a
tautological family sheaves $\cF$ on $S$ parametrized by $\mx$,
i.e.~such that $\cF|_{S\tm[F]}\cong F$. We state below a simple condition
ensuring that $\mx$ is a fine moduli space:  the verification
that semistability implies stability is left to the reader, the
existence of a tautological sheaf follows from~\cite{Ma (6.11),Mk2 (A.7)}.
\proclaim{(0.3) Criterion}
Assume that for $[F]\in \mx$
$$\gcd\left\{ \rk(F), c_1(F)\cdot H, \chi(F) \right\}=1.$$
Then $\mx$ is a fine moduli space. 
\endproclaim
\head
1. Outline of the proof of Theorem~(0.1).
\endhead
The moduli space $\cM$ appearing in Theorem~(0.1) parametrizes sheaves
with fixed rank $r$, $c_1\in H^{1,1}(S;\ZZ)$, and $c_2\in H^4(S;\ZZ)$. Let
$\xi$ be a set of sheaf data with $\rk(\xi)=r$, $c_1(\xi)=c_1$,
$c_2(\xi)=c_2$. Since $\cM$ is a locally-trivial fibration over
$Pic^{c_1}(S)$, with fiber isomorphic to $\mx$, Theorem~(0.1) is equivalent
to the analogous statement obtained replacing $\cM$ by $\mx$. (Of course
the expected dimension of $\mx$ is obtained subtracting $h^1(\cO_S)$ from
the expected dimension of $\cM$.) We will outline the proof of the
statement for $\mx$: hence from now on we will only deal with $\mx$, the
moduli space with fixed determinant.   
\bsk

\n
{\bf Deformation theory and twisted endomorphisms.}
\hskip 2mm
References for deformation theory are~\cite{A,F2,Mk1,ST}. Let
$[F]\in\mx^{st}$, i.e.~$F$ is stable. The germ of $\mx$ at $F$ is isomorphic
to $\Def^0(F)$, the universal deformation space of $F$ "with fixed
determinant" (i.e.~it classifies deformations of $F$ which do not change the
isomorphism class of $\det F$). To describe $\Def^0(F)$ we need the
traceless $\Ext$-groups. If $L$ is a line-bundle on $S$ we set 
$$\Ext^q(F,F\ot L)^0:=
\ker\left(\Ext^q(F,F\ot L)\brel \Tr\over\lra H^q(L)\right).$$
The trace $\Tr$ is defined in~\cite{DL}; if $F$ is locally-free then  
$$\Ext^q(F,F\ot L)^0=H^q(End_0(F)\ot L),$$
 where $End_0(F)$ is the sheaf of
traceless endomorphisms of $F$. We set 
$$h^q(F,F\ot L)^0:=\dim \Ext^q(F,F\ot L)^0.$$

The tangent space to $\Def^0(F)$ is canonically identified with
$\Ext^1(F,F)^0$. There is a Kuranishi map
$$\Ext^1(F,F)^0\supset U\brel\Phi\over\lra\Ext^2(F,F)^0,$$
defined on  an open neighborhood $U$ of the origin,  such that $\Def^0(F)$ is
the germ at the origin of $\Phi^{-1}(0)$. Thus 
$$\eqalign{\dim_{[F]}\mx\ge & \dim\Ext^1(F,F)^0-\dim\Ext^2(F,F)^0
=\chi(F,F)^0=\cr
&2\rk(\xi)\D_{\xi}-\left(\rk(\xi)^2-1\right)\chi(\cO_S)
=:\text{exp.dim.}\left(\mx\right).\cr}$$
In fact the first equality holds because since $F$ is stable
$\Hom(F,F)^0=0$, and the second eqality is just Riemann-Roch. The
obstruction space $\Ext^2(F,F)^0$ is Serre dual to $\Hom(F,F\ot K)^0$, hence
we have the following. 
\proclaim{Criterion}
Assume the locus of $[F]\in \mx^{st}$ such that
$$h^0(F,F\ot K)^0>0$$
has dimension strictly smaller than the expected dimension of $\mx$. Then
$\mx^{st}$ is a reduced local complete intersection scheme of dimension
the expected one.  
\endproclaim
For $L$ a line-bundle on $S$, let 
$$W_{\xi}^L:=\{[F]\in\mx^{st}|\ h^0(F,F\ot L)^0>0\}.$$
Theorem~(0.1), except for the statement about irreducibility, follows
essentially from the following result.
\proclaim{(1.1) Theorem~\cite{O2}} 
There exist numbers $\l_0'(\rk(\xi),S,H,L)$, $\l_1(\rk(\xi),S,H)$ and
$\l_2(\rk(\xi))$, with $\l_2(\rk(\xi))<2\rk(\xi)$, such that
$$\dim W_{\xi}^L\le \l_2\D_{\xi}+\l_1\sqrt{\D_{\xi}}+\l_0'.$$
\endproclaim
Indeed the theorem implies  $\dim W_{\xi}^K$ is
strictly less than the expected dimension,  if $\D_{\xi}$ is large enough: 
by the previous criterion $\mx^{st}$ is reduced of the expected
dimension. To deal with $\left(\mx-\mx^{st}\right)$, i.e.~strictly
semistable sheaves, one needs some dimension counts: this is a
technical point. For simplicity we will usually ignore strictly semistable
sheaves: as a first approximation the reader may assume $\mx$ is
a fine moduli space (see~(0.3)). Similarly the statement in Theorem~(0.1)
that the locus parametrizing slope-stable vector-bundles is dense follows
from Theorem~(1.1) together with  a result of Jun Li~\cite{L1, Appendix}. 
\remark{Remark}
The coefficients in the above theorem can be computed
explicitly: they depend on $(S,H,L)$ only via intersection numbers, in
particular they are constant for families of polarized surfaces.
In~~\cite{O2} there are some explicit lower bounds for $\dx$ ensuring
$\mx$ is reduced of the expected dimension. Donaldson~\cite{D,F2,Z} proved 
Theorem~(1.1) for rank two: his coefficient of $\dx$ is $3$, which is better
than our $\l_2(2)=23/6$, but the other coefficients are not explicit. 
We will see later (see~(2.6)) how to use Theorem~(1.1) with
choices of $L$ different from $K_S$. 
\endremark
In  this section we will sketch a proof of Theorem~(1.1) and we will give
the argument for proving (asymptotic) irreducibility. 
\subhead
The boundary
\endsubhead
If $X\ss\mx$, {\it the boundary} $\del X$ consists of the subset of points
parametrizing singular (i.e.~not locally-free) sheaves. Our approach to the
proof of Theorem~(1.1) is to show that any closed subset of $\mx$ of
relatively small codimension has non-empty  boundary. More
precisely we prove the following result.
\proclaim{(1.2) Theorem}
 There exists $\l_0\left(\rk(\xi),S,H\right)$ such that if $X$ is a closed
irreducible subset of $\mx$ with
$$\dim X> \l_2\dx+\l_1\sqrt{\dx}+\l_0,$$
then $\del X$ is non-empty. (Here $\l_2$, $\l_1$ are as in 
Theorem~(1.1).)
\endproclaim
We will illustrate the implication
$\text{Theorem(1.2)} \Longrightarrow \text{Theorem(1.1)}$
by proving the following.
\proclaim{(1.3) Proposition}
Assume Theorem~(1.2) holds. Let $r\ge 2$ be an integer and $D$ be a
divisor on $S$. Suppose the following:  if a torsion-free sheaf $F$ with
$\rk(F)=r$ and $\det F\cong\cO_S(D)$ is semistable then it is slope-stable
(e.g.~if $D\cdot H$ and $r$ are coprime). If $L$ is a line-bundle on
$S$ then 
$$\dim W_{\xi}^L<
\text{exp.dim.}\left(\mx\right)=2r\dx-(r^2-1)\chi(\cO_S)$$
for all sheaf data $\xi$ such that $\rk(\xi)=r$, $\det(\xi)\cong\cO_S(D)$,
and $\dx>>0$. 
\endproclaim
Before proving the above proposition we need some preliminaries on
double-duals. Let $F$ be a torsion-free sheaf on $S$. Since $\dim S=2$ and
$S$ is smooth the double-dual $F^{**}$ is locally-free~\cite{OSS}, since $F$
is torsion-free the natural map $F\to F^{**}$ is an injection. Thus we get a
canonical exact sequence 
$$0\to F\to F^{**}\to Q(F) \to 0.$$
The   lenght $\ell(Q(F))$ is finite. We have
$$\rk(F^{**})=\rk(F),\qquad \det(F^{**})=\det(F),
\qquad c_2(F^{**})=c_2(F)-\ell(Q(F)).\tag 1-4$$
In particular $F$ is slope-stable if and only if so is $F^{**}$. Now let $\xi$
be a set of sheaf data as in the statement of Proposition~(1.3). 
Let $X\ss\mx$ be a closed irreducible subset. If $[F]\in\del X$ then
by our hypothesis $F^{**}$ is slope-stable. Thus $[F^{**}]\in\cM_{\xi'}$,
where $\xi'$ is determined by~(1.4). The double-duals $F^{**}$, for $[F]$
varying in $\del X$, are not parametrized by a single moduli space: in
general $c_2(F^{**})$ will vary with $[F]$. However   $\del X$ is stratified
by the {\it double-dual} strata:  if $[F]$ varies in a single stratum then
$[F^{**}]$ varies (algebraically) in a single moduli space (each stratum is
locally closed). Let $Y\ss \del X$ be an irreducible component of the open
stratum: we set  
$$Y^{**}:=\{[F^{**}]|\ [F]\in Y\},\qquad\ell:=c_2(\xi)-c_2(\xi').$$ 
We will need an inequality between the dimensions of $X$ and  $Y^{**}$.
First of all, considering short locally-free resolutions of sheaves
parametrized by $X$, one gets that
$$\cod(\del X,X)\le r-1.$$
Secondly, a sheaf parametrized by $Y$ is determined by the isomorphism
class of its double-dual, i.e.~a point $[E]\in Y^{**}$, plus the choice of a
quotient $E\to Q$, where $\ell(Q)=\ell$. A theorem of Jun Li~\cite{L1,
Appendix} asserts that the generic such quotient is isomorphic to
$\bigoplus_{i=1}^{\ell}\CC_{P_i}$. Putting together these facts one obtains
the following.
\proclaim{Lemma}
Keeping notation as above,
$$\dim Y^{**}=\dim X-2r\ell+(r-1)(\ell-1)+\e,\tag 1-5$$
where $\e\ge 0$. If $\e=0$ then $\del X$ contains the isomorphism class
of all sheaves $F$ fitting into an exact sequence
$$0\to F\to E\brel\phi\over\lra\bigoplus_{i=1}^{\ell}\CC_{P_i}\to 0,$$
where the point $[E]\in Y^{**}$, the points $P_i\in S$, and the
surjection $\phi$ are chosen arbitrarily.  
\endproclaim
The reader should notice that $\text{exp.dim.}\left(\cM_{\xi'}\right)= 
\left(\text{exp.dim.}\left(\cM_{\xi}\right)-2r\ell\right)$.
The proof of the proposition will go roughly as follows. Starting from
$X_0:=W_{\xi}^L$ we will repeatedly apply Theorem~(1.2) and
construct as above $Y_0\ss\del X_0$, $X_1:=Y_0^{**}$, $Y_1\ss\del
\ov{X}_1$, and so on. We will show that in most cases the quantity $\e$
of Inequality~(1.5) is strictly positive. This progressively "inflates" the
dimension of $X_i$, until it becomes too big, giving a contradiction. We
still have to introduce a key ingredient in this argument, namely an a priori
bound on the amount by which the actual dimension of a moduli space can
exceed the expected dimension. This follows from a bound for the
number of sections of semistable sheaves, obtained by Simpson~\cite{S,
Cor.~(1.7)}. For the purposes of this proof we only need to know that there
exists $e_L(r,S,H)$ such that 
$$h^0(F,F\ot L)^0\le  e_L(r,S,H)\tag 1-6$$
for all slope-semistable sheaves $F$ with $\rk(F)=r$; the point is that
$e_L$ is independent of the discriminant $\D_F$. Under the hypotheses
of~(1.3) we have $\mx=\mx^{st}$, hence deformation theory gives
$$\dim\mx\le\text{exp.dim.}(\mx)+e_K.\tag1-7$$
\demo{Proof of Proposition~(1.3)}
Let $\D_0$ be so large that  
$$\text{exp.dim.}(\mx)>\l_2(r)\dx+\l_1(r,S,H)\sqrt{\dx}+\l_0(r,S,H)$$
for all $\xi$ with $\dx\ge\D_0$ (here $\l_2$, $\l_1$, $\l_0$ are as in
Theorem~(1.2)). By Theorem~(1.2) we have the following.
\proclaim{(1.8)}
Assume $\dx\ge\D_0$. If $X$ is a closed irreducible subset of $\mx$, 
with $\dim X\ge\text{exp.dim.}(\mx)$, then $\del X\not=\es$.
\endproclaim
Now assume   
$$\dx>\D_0+e_L+e_K.\tag 1-9$$
Let's show that $\dim W_{\xi}^L<\text{exp.dim.}(\mx)$. Suppose the
contrary, and let $X_0\ss W_{\xi}^L$ be an irreducible component  with
$\dim X_0\ge\text{exp.dim.}(\mx)$. By~(1.9) and~(1.8) $\del X_0\not=\es$.
Let $Y_0\ss\del X_0$ be an irreducible component of the open double-dual
stratum, and set $X_1:=Y_0^{**}$. If $\del\ov{X}_1\not=\es$ ($\ov{X}_1$ is
the closure of $X_1$ in the appropriate moduli space) we repeat the process,
i.e.~we consider $Y_1\ss\del\ov{X}_1$, $X_2:=Y_1^{**}$, and continue until
we reach $X_n$ such that $\del\ov{X}_n=\es$. By Formula~(1.5) we have 
$$\dim X_{i+1}=\dim X_i-2r\ell_i+(r-1)(\ell_i-1)+\e_i,$$
with the obvious notation. Let $\cM_{\xi_n}$ be the moduli space to which
$X_n$ belongs. The formula above gives that 
$$\dim X_n=\dim X_0-2r\sum_{i=0}^{n-1}\ell_i+
(r-1)\left(\sum_{i=0}^{n-1}(\ell_i-1)\right)+\sum_{i=0}^{n-1}\e_i\ge
\text{exp.dim.}\left(\cM_{\xi_n}\right).\tag{$*$}$$
In fact the sum of the first two terms equals
$\text{exp.dim.}\left(\cM_{\xi_n}\right)$, and the remaining terms are
non-negative. Since we are assuming $\del \ov{X}_n=\es$,
we conclude by~(1.8) that $\D_{\xi_n}<\D_0$. Since 
$\dx-\D_{\xi_n}=\sum_{i=0}^{n-1}\ell_i$,
$$\dx-\D_0\le\sum_{i=0}^{n-1}\ell_i.$$
Manipulating the second term of~($*$), and applying the above inequality we
get 
$$\dim X_n\ge\text{exp.dim.}\left(\cM_{\xi_n}\right)+
\dx-\D_0+\sum_{i=0}^{n-1}(\e_i-1).\tag 1-10$$
Now comes the key observation.
\proclaim{Claim}
Let $h_L(X_i):=\min\{h^0(F,F\ot L)^0|\ [F]\in X_i\}$. Then:
\roster
\item $0<h_L(X_i)\le h_L(X_{i+1})$ for $i=0,\ldots,n-1$.
\item If $\e_i=0$ then $h_L(X_i)<h_L(X_{i+1})$.
\endroster
\endproclaim
\demo{Proof of the claim}
To prove Item~(1) it suffices to show that $h_L(X_i)\le h_L(X_{i+1})$ for all
$i$, because $h_L(X_0)>0$ by definition. If $F$ is a torsion-free
sheaf on $S$ there is a canonical injection 
$$\rho\cl \Hom(F,F\ot L)\hra \Hom(F^{**},F^{**}\ot L)$$
which commutes with the trace, hence it defines also an injection of the
traceless $\Hom$ groups. As is easily seen this implies~(1): Indeed let
$[F]\in Y_i\ss\del X_i$ be a generic point; by upper-semicontinuity
$h^0(F,F\ot L)^0\ge h_L(X_i)$, hence  $h^0(F^{**},F^{**}\ot L)^0\ge h_L(X_i)$.
Since $F^{**}$ is a generic point of $X_{i+1}$ we have 
$h_L(X_{i+1})=h^0(F^{**},F^{**}\ot L)^0$; we have proved Item~(1). 
Now let's prove Item~(2). The hypothesis together with
Equation~(1.5) implies that $\del X$ contains all sheaves $F$
fitting into an exact sequence
$$0\to F\to E\brel\phi\over\to\bigoplus_{j=1}^{\ell_i}\CC_{P_j}\to 0,$$
where $[E]\in X_{i+1}$. Clearly $E=F^{**}$, thus the map $\rho$ realizes
$\Hom(F,F\ot L)^0$ as a subgroup of $\Hom(E,E\ot L)^0$; an element $f\in 
\Hom(E,E\ot L)^0$ belongs to the image of $\rho$ if and only if
$$\text{$f(\Ker\phi_j)\ss\Ker\phi_j\ot L$, for
$j=1,\ldots,\ell_i$,}\tag{$\bu$}$$ 
 where $\phi_j$ is the restriction of $\phi$ to the fiber over $P_j$. 
Now let $[E]\in X_{i+1}$ be generic: by
upper-semicontinuity $h^0(E,E\ot L)^0=h_L(X_{i+1})$, and the latter is
non-zero by Item~(1). Let $f\in\Hom(E,E\ot L)^0$ be non-zero:
since $f$ is not a scalar endomorphism at the generic point, we can choose
$\phi$ (in fact the generic $\phi$ will do) so that~($\bu$) does not hold,
i.e.~$\rho$ is not surjective. Hence for generic $[F]\in Y_i\ss\del X_i$ we
have $h^0(F,F\ot L)<h_L(X_{i+1})$. By upper-semicontinuity
$h_L(X_i)<h_L(X_{i+1})$. 
\qed 
\enddemo
Let's conclude the proof of Proposition~(1.3). Since $h_L(X_i)\le e_L$ for
all $i$ by Simpson's bound~(1.6), the claim implies that
$\sum_{i=0}^{n-1}(\e_i-1)\ge-(e_L-1)$. By~(1.10) we conclude that 
$$\dim X_n\ge\text{exp.dim.}\cM_{\xi_n}+\dx-\D_0-e_L+1.$$
Since $\dx$ satisfies~(1.9) this inequality contradicts~(1.7). 
\qed
\enddemo
We will  sketch a proof of Theorem~(1.2). First we need to
discuss determinant bundles on the moduli space.
\subhead
Determinant bundles
\endsubhead
References for this section are~\cite{LP,L1,FM (5.3.2)}. Assume $\mx$ is
fine, thus there is a tautological sheaf $\cF$ on  $S\tm\mx$. Let $C\ss S$
be a smooth irreducible curve. Choose a vector-bundle $A$ on $C$ with the
property that 
$$\chi\left(F|_C \ot A\right)=0
\qquad \text{for all $[F]\in\mx$.} \tag 1-11$$ 
The restriction  $\cF|_{C\tm\mx}$ is flat over $\mx$, hence by the
theory of determinant  line-bundles~\cite{KM} it makes sense to set
$$\cL(\cF,C,A):=det R q_!\left(\cF\ot p^*A\right)^{-1},$$
where $p,q$ are the projections of $C\tm\mx$ to $C$ and $\mx$
respectively. Since $A$ satisfies~(1.11) the determinant line-bundle is
independent of the choice of a tautological sheaf. There is a natural section
of $\cL(\cF,C,A)$ whose zero-locus is supported on the subset
parametrizing sheaves $F$ such that $h^0(F|_{C}\ot A)>0$ (of course it
might  be that this section vanishes identically on $\mx$).
Applying the Grothendieck-Riemann-Roch Theorem (for Chow groups), one
gets the equality 
$$c_1\left(\cL(\cF,C,A)\right)=\rk(A)\pi_*\left[\left(c_2(\cF)-
{\rk(\cF)-1\over 2\rk(\cF)}c_1(\cF)^2\right)\cdot C\right],\tag 1-12$$
where $\pi\cl S\tm\mx\to\mx$ is the projection (use Equation~(1.11)).
The above formula shows that the isomorphism class of $\cL(\cF,C,A)$ only
depends on the linear equivalence class $[C]$. Furthermore, since the
right-hand side of~(1.12) is linear in $[C]$, we can define $\cL(\cF,[C],A)$
for an arbitrary divisor class $[C]$ . To get rid of the dependence
from $\rk(A)$ we set 
$$\cL([C]):={1\over \rk(A)}\cL(\cF,[C],A).$$
Thus we get a well-defined map $\cL\cl \Pic(S)\to\Pic(\mx)\ot\QQ$. We
set $\cL(n):=\cL([nH])$; as is easily verified $\cL(n)$ is a
line-bundle for all $n$ divisible by $\rk(\cF)$. 
For simplicity we have assumed that $\mx$ is fine, but in fact the map $\cL$
can be defined without this assumption~\cite{L1,LP}: the domain of $\cL$
will be a certain subspace of $\Pic(S)$ which always includes $\ZZ[H]$. 
Historically Donaldson~\cite{D} was the first to study the determinant
line-bundle: his goal was to prove  that the  polynomial invariants of 
algebraic surfaces are non zero. The following theorem gives an important
property of  $\cL(n)$~\cite{LP,L1}.
\proclaim{(1.13) Theorem (Le Potier - J.~Li)}
Let $n$ be  sufficiently large and divisible by $\rk(\xi)$ (in particular
$\cL(n)$ is line-bundle). Then the complete linear system $|\cL(n)|$ is
base-point free, and it defines an embedding of the subset of $\mx$
parametrizing $\mu$-stable  locally-free sheaves.
\endproclaim
We will use the following.
\proclaim{(1.14) Corollary}
Let $X\ss \mx$ be a closed irreducible subset. If the generic point of $X$
parametrizes a $\mu$-stable locally-free sheaf then 
$$c_1(\cL(n))^{\dim X}\cdot X >0.$$
\endproclaim
The rational line-bundle $\cL(n)$ is related to the theta-divisor on the
moduli space of vector-bundles on $C\in |nH|$, as follows. Let
$\cA_C\ss\mx$ be the subset parametrizing sheaves whose restriction to
$C$ is locally-free and stable; restriction  defines a morphism
$$\rho\cl \cA_C\to \cU(C;\rk(\xi),\det(\xi)|_C)$$
to the moduli space of semistable vector-bundles on $C$ (with fixed
determinant). If $\T$ is the theta-divisor on $\cU(C;\rk(\xi),\det(\xi)|_C)$,
then
$$\rho^*\T\sim \l c_1\left(\cL(n)\right),\tag 1-15$$
where $\l$ is a positive integer. 
\subhead
The proof of Theorem~(1.2)
\endsubhead
For simplicity we assume $\mx$ is a fine moduli space. To
lighten notation we set $r=\rk(\xi)$ and $\D=\dx$. The proof is by
contradiction. So let's assume $X\ss\mx$ is an irreducible closed subset
with $\del X=\es$. If $C\ss S$ is an irreducible smooth curve we set
$$X_C:=\{[F]\in X|\ \text{$F|_C$ is not stable}\}.$$
A key observation is that under certain hypotheses $X_C$ is non-empty
\proclaim{(1.16) Proposition}
Keep notation and hypotheses as above (in particular $\del X=\es$). Suppose
$n$ is a positive integer such that
$${r^2-1\over 2}H^2n^2+{r^2-1\over 2}K\cdot Hn<\dim X.\tag 1-17$$
If $C\in |nH|$ is a smooth curve, then $X_C$ is non-empty, and moreover
$$\dim X_C\ge 
\dim X-{r^2\over 8}H^2n^2-{r^2\over 8}K\cdot Hn-{r^2\over 4}.
\tag 1-18$$ 
\endproclaim
\demo{Proof}
Assume that $X_C=\es$. Then associating to $[F]\in X$ the S-equivalence
class of $F|_C$ we get a well-defined morphism
$$\rho\cl X\to \cM(C;\xi),$$
where $\cM(C;\xi)$ is the moduli space of rank-$r$ semistable
vector-bundles on $C$ with determinant $\det(\xi)|_C$. Since the left-hand
side of~(1.17) equals $\dim\cM(C;\xi)$, we have
$$\left(\rho^*\Theta\right)^{\dim X}=0,$$
where $\Theta$ is the theta-divisor. By Equation~(1.15) and
Corollary~(1.14) we conclude that the generic point (hence all points) of
$X$ parametrizes a sheaf which is not slope-stable. This contradicts our
assumption that $X_C=\es$: in fact it follows directly from the definition
of slope-stability that if $F|_C$ is stable (where $C\in |nH|$), then $F$ is
slope-stable. This proves $X_C\not=\es$. Once we know $X_C\not=\es$, 
Inequality~(1.18) follows from a straightforward dimension count.
\qed
\enddemo
Now assume we are in the situation of Proposition~(1.16). Choose $[F]\in
X_C$, and let
$$0\to\cL_0\to F|_C\brel g\over \to \cQ_0\to 0\tag 1-19$$
be a destabilizing sequence for $F|_C$ (with $\cQ_0$ locally-free). Let $E$
be the locally-free sheaf on $S$ defined by the following exact sequence
(an {\it elementary modification}) 
$$0\to E\to F\brel g\over\to\i_*\cQ_0\to 0,$$
where $\i\cl C\hra S$ is the inclusion. Restricting to $C$ the above
sequence we get an exact sequence
$$0\to \cQ_0\ot\cO_C(-C)\to E|_C\brel f_0\over \to \cL_0\to 0.$$
Let $Y_F:=Quot(E|_C;\cL_0)$ be the Quot-scheme parametrizing quotients of
$E|_C$ with Hilbert polynomial equal to that of $\cL_0$. For $y\in Y_F$ we
let $\cG_y$ be the torsion-free sheaf on $S$ defined by the elementary
modification
$$0\to\cG_y\to E\brel f_y\over\to \i_*\cL_y\to 0,$$
where $f_y$ is given by the quotient of $E|_C$ parametrized by $y$. The
sheaves $\cG_y$ fit into a family parametrized by $Y_F$. One
easily verifies that:
\roster
\item
There is a natural isomorphism $\cG_0\ot\cO_S(C)\cong F$.
\item 
The sheaf $\cG_y$ is singular if and only if so is $\cL_y$
(i.e.~if $\cL_y$ has torsion). 
\endroster
Let's assume for the moment that $\cG_y$ is stable for all $y\in Y_F$.
Then, setting $\cF_y:=\cG_y\ot\cO_S(C)$,   the
family $\{\cF_y\}$ defines a classifying morphism
$$\vf\cl Y_F\to \mx,$$
and by Item~(1) we have $\vf(0)=[F]\in X$. We will arrive at a contradiction
if we show that there exists $y\in\vf^{-1}X$ such that $\cL_y$ is singular;
indeed this implies $\cG_y$ is singular by Item~(1), hence $\cF_y$ is also
singular, and  thus $\vf(y)\in\del X$, contradicting the assumption $\del
X=\es$. The following elementary  result is proved~\cite{O2}. 
\proclaim{(1.20) Lemma}
Let $\S\ss Y_F$ be a closed irreducible subset with $\dim\S>r^2/4$.
There exists $y\in\S$ such that $\cL_y$ is singular. 
\endproclaim
To apply the lemma we notice that
$$\dim\vf^{-1}X\ge \dim_0 Y_F+\dim X-\dim T_{[F]}\mx.\tag 1-21$$
For the dimension of the Quot-scheme $Y_F$ we have
$$\align
\dim_0Y_F &\ge 
\chi\left(Hom(\cQ_0\ot\cO_C(-C),\cL_0)\right)\\
& =\rk(\cL_0)\rk(\cQ_0)\left(\mu(\cL_0)-\mu(\cQ_0)+C^2+1-g(C)\right)\\
& \ge \rk(\cL_0)\rk(\cQ_0)\left({1\over 2}C^2-{1\over 2}C\cdot K\right)\\
& \ge {1\over 2}(r-1)\left(H^2n^2-K\cdot Hn\right).
\endalign$$ 
(The second inequality holds because~(1.19) is a destabilizing sequence.) 
Feeding the inequality for $\dim Y_F$ together with~(1.7) into~(1.21),
and applying Lemma~(1.20) we get the following.
\proclaim{(1.22)}
Assume $\dim X$ satisfies~(1.17). Assume also that $\cG_y$ is stable for
all $y\in Y_F$. If
$$\dim X> 2r\D-(r^2-1)\chi(\cO_S)+e_K-{1\over 2}(r-1)
\left(H^2n^2-K\cdot Hn\right),\tag 1-23$$
then there exists $y\in \vf^{-1}X$ with $\cL_y$ singular, and hence $\del
X\not=\es$. 
\endproclaim
To deal with the condition that $\cG_y$ be stable for all $y\in Y_F$ we
want to choose $[F]\in X_C$ which is "very stable", i.e.~such that for all
subsheaves $E\ss F$ with $\rk(E)<\rk(F)$, 
$$\mu(E)<\mu(F)-C\cdot H=\mu(F)-H^2n.$$
Carrying out some dimension counts and using~(1.18) one shows  it
suffices that 
$$\dim X-{r^2\over 8}H^2n^2-{r^2\over 8}K\cdot Hn-{r^2\over 4}>
(2r-1)\D+(2r-1)(r-1)^2H^2n^2+O(n).\tag 1-24$$
(For this we must assume $|H|$ is base-point free.) If $r=2$ a weaker
inequality is required~\cite{O2}. At this point we have all the elements
needed to prove Theorem~(1.2). If we can find $n$ such that
Inequalities~(1.17)-(1.23)-(1.24) hold, then the argument sketched above
shows that $\del X\not=\es$. It is an easy exercise to determine a lower
bound on $\dim X$ guaranteeing such $n$ exists. The reader can check that
the coefficient of $\D$ can be taken to be
$$\l_2(r)=2r-{4(r-1)\over 16r^3-39r^2+36r-12}.$$
(If $r=2$ one can improve the estimates and get $23/6$ rather
than $31/8$.)   

The lower bound on $\dx$ ensuring that $\mx$ is reduced of the expected
dimension can be computed explicitly. This has been carried out in~\cite{O2}
for $\rk(\xi)=2$, when $K$ is ample and $H=K$. The lower bound is of the
form $(cost.)K^2$. One can ask for sharp bounds:
\remark{Question}
Assume $S$ is minimal of general type. Is $\mx$ reduced of the expected
dimension when
$$\dx>\rk(\xi)(p_g+1),$$
for polarizations sufficiently close to $K$ ?
\endremark
Notice that we must restrict the choice of polarization $H$ or else the
answer is certainly negative~(see~\cite{O3 (5b.24)}): sufficiently close
means that for a sheaf with Chern classes defined by $\xi$ slope-stability
(instability) for $H$ and $K$ coincide.
\subhead
Irreducibility
\endsubhead
We give the  argument of Gieseker and Li~\cite{GL1} which proves that
$\mx$ is irreducible for large enough $\dx$; we will make some
simplifying assumptions (as in Proposition~(1.3)) in order to avoid some
minor technical problems.
\proclaim{(1.25) Theorem}
Let $r\ge 2$ be an integer, and $D$ be a divisor on $S$. Suppose  that
every rank-$r$ torsion-free semistable sheaf $F$ on $S$ with $\det
F\cong\cO_S(D)$ is actually slope-stable (e.g.~if $D\cdot H$ and $r$ are
coprime). There exists $\D_1$ such that if $\xi$ is a set of sheaf data
with 
$$\rk(\xi)=r,\qquad \det(\xi)\cong\cO_S(D),\qquad \dx>\D_1,$$
then $\mx$ is irreducible.
\endproclaim
This section is devoted to proving Theorem~(1.25). We will  always
assume that $\xi$ is a set of sheaf data satisfying the hypotheses of the
theorem. 

Let $\xi_0$ be a set of sheaf data
(with $\rk(\xi_0)=r$, $\det(\xi_0)=\cO_S(D)$),  and let
$X_1,\ldots,X_n$ be the irreducible components of $\cM_{\xi_0}$.
For $\ell$ a positive integer, and $i=1,\ldots,n$, we let $Y_i^{\ell}$ be the
locus of moduli (in the appropriate moduli space) of sheaves $F$ fitting into
an exact sequence 
$$0\to F\to E\to\op_{j=1}^\ell\CC_{P_j}\to 0,$$
where $[E]\in X_i$ is an arbitrary point with $E$ locally-free, and the
$P_j$'s are pairwise distinct. 
\proclaim{Lemma}
Keep notation as above. There exists $\D_{\xi_0}$ such that  the
following holds. If $\dx>\D_{\xi_0}$,  and $\ell:=c_2(\xi)-c_2(\xi_0)$, then
any irreducible component of $\mx$ contains one (at least) of the
$Y_i^{\ell}$. Furthermore $\mx$ is smooth at the generic point of each of the
$Y_i^{\ell}$. 
\endproclaim
\n
\demo{Sketch of proof}
Let $\D_0$ be as in~(1.8): hence if $\dx\ge \D_0$ all irreducible components
of $\mx$ have non-empty boundary.  Increasing $\D_0$ if necessary, we can
assume by Proposition~(1.3) that moduli spaces $\mx$ with $\dx\ge \D_0$
are reduced of the expected dimension. A simple application of
Inequality~(1.5) will show that if $\D_{\xi_0}$ is sufficiently larger than
$\D_0$ the following holds. Assume $\dx>\D_{\xi_0}$, and let $V$ be any
irreducible component of $\mx$. Then there exists an irreducible component
$V'$ of $\cM_{\xi'}$, where 
$$\rk(\xi')=r,\qquad \det(\xi')=\cO_S(D),\qquad c_2(\xi')=c_2(\xi)-1,$$
such that $V$ contains the moduli point of any sheaf $F$ fitting into an
exact sequence
$$0\to F\to E\to \CC_P\to 0,$$
where $[E]$ is an arbitrary point of $V'$ with $E$ locally-free. Applying
this same result to $\cM_{\xi'}$ and the irreducible component $V'$, and
so on  all the way down to $\D_{\xi_0}$, one gets the first statement of the
lemma.  The second statement holds because $\D_{\xi_0}\ge\D_0$, and
hence the generic point $[E]$ of any irreducible component of $\cM_{\xi_0}$
has vanishing obstruction space (i.e.~$H^0(End_0(E_i)\ot K)=0$), and hence
so does any sheaf whose double-dual is isomorphic to $E$.  
\qed
\enddemo
Fix  $\xi_0$ as in the above lemma; then for
$i=1,\ldots,n$ there is only one irreducible component of $\mx$ containing
$Y_i^{\ell}$, and since each component contains at least one $Y_i^{\ell}$, 
$$\#\text{irr.comp.}\left(\mx\right)\le 
\#\text{irr.comp.}\left(\cM_{\xi_0}\right).$$
We will prove Theorem~(1.25) by showing that if $\ell\gg 0$ then
the $Y_i^{\ell}$ all belong to the same irreducible component. Choose
$[E_i]\in X_i$, for $i=1,\ldots,n$, with $E_i$ locally-free and with
vanishing obstruction space . Thus 
if $[F_i]\in Y_i^{\ell}$ lies over $[E_i]$, i.e.~$F_i^{**}\cong E_i$, the
moduli space $\mx$  is smooth at
$[F_i]$, in particular the unique  irreducible component containing all of
$Y_i^{\ell}$ must contain any irreducible subset through $[F_i]$. We will
construct (for $\ell\gg 0$) an irreducible subset $W\ss\mx$ containing
$[F_1],\ldots,[F_n]$; thus $\mx$ must be
irreducible. The subset $W$ is defined as follows. Let $n$ be an integer
such that $E_i\ot\cO_S(n)$ is generated by global sections, for
$i=1,\ldots,n$. Choosing $(r-1)$ generic sections of $E_i\ot\cO_S(n)$ we
see that $E_i$ fits into an exact sequence  
$$0\to \cO_S(-nH)^{(r-1)}\to E_i\to I_{Z_i}(D+(r-1)nH)\to 0,$$
where $Z_i$ is some zero-dimensional subscheme of $S$. Choosing
appropriately the surjection 
$$E_i\to\op_{j=1}^{\ell}\CC_{P^i_j}$$
 whose kernel is $F_i$, we see that $Y_i$ contains the moduli point of a
sheaf $F_i$ fitting into an exact sequence
$$0\to\cO_S(-nH)^{(r-1)}\to F_i\to I_{\wt{Z}_i}(D+(r-1)nH)\to 0,$$
where $\wt{Z}_i=Z_i\cup\{P^i_1,\ldots,P^i_{\ell}\}$. Thus $F_i$ corresponds
to a non-zero class in
$$\Ext^1\left(I_{\wt{Z}_i}(D+(r-1)nH),\cO_S(-nH)^{(r-1)}\right).$$
If $\ell$ is large enough and the points $P^i_1,\ldots,P^i_{\ell}$ are generic,
$$\multline
\dim\Ext^1\left(I_{\wt{Z}_i}(D+(r-1)nH),\cO_S(-nH)^{(r-1)}\right)\\
=-\chi\left(I_{\wt{Z}_i}(D+(r-1)nH),\cO_S(-nH)^{(r-1)}\right).
\endmultline\tag{$\dag$}$$
Let $d=\ell(\wt{Z}_1)=\ldots=\ell(\wt{Z}_n)$, and let $U\ss S^{[d]}$ be the
open subset of the Hilbert scheme parametrizing length-$d$ subschemes
$\wt{Z}$ of $S$ such that~($\dag$) holds with $\wt{Z}_i$ replaced by
$\wt{Z}$. We define $W\ss\mx$ to be the subset parametrizing sheaves $F$
which fit into an exact sequence  
$$0\to\cO_S(-nH)^{(r-1)}\to F\to I_{\wt{Z}}(D+(r-1)nH)\to 0,$$
for some $\wt{Z}\in U$. By construction $[F_i]\in W$ for $i=1,\ldots,n$. 
Since $W$ is an open subset of a bundle of projective
spaces over $U$, it is irreducible. This finishes the proof of
Theorem~(1.25). 
\remark{Remark}
Notice that all the steps of the above proof can easily be made effective,
except for the choice of $n$ such that $E_i\ot\cO_S(nH)$ is generated by
global sections. In~\cite{O2} there are some effective results for complete
intersections.   
\endremark
\head
2.  Two-forms on the moduli space.
\endhead
Let $B$ be a smooth variety, and $\cF$ be a family of torsion-free sheaves
on $S$ parametrized by $B$. For $b\in B$ we set $S_b:=S\tm\{b\}$ and
$\cF_b:=\cF|_{S_b}$.  We assume the isomorphism class of $\det\cF_b$ is
independent of $b$. Given a two-form $\o\in\G(\O^2_S)$ we will define a
two-form $\o_{\cF}\in\G(\O^2_B)$. 
First  recall~\cite{Mm1} that given a codimension-two cycle $\cZ\in
Z^2(S\tm B)$   transverse to the projection $q\cl S\tm B\to B$
(i.e.~$\cZ=\sum_in_i\cZ_i$ where each $\cZ_i$ is a subvariety
intersecting the generic $S_b$ in a finite set of points)
we can associate to it a two-form $\o_{\cZ}$ on $B$. Explicitely,  let
$q_i\cl\cZ_i\to B$ be the restriction of $q$, and $p\cl S\tm B\to S$ be the
projection, then  
$$\o_{\cZ}:=\sum_in_iq_{i,*}(p^*\o|_{\cZ_i}).\tag 2-1$$
Some care must be taken in defining the push-forward at points $b\in B$
over which $q_i$ is not \'etale: we can circumvent this problem by
considering the universal case, i.e.~$B=S^{[d]}$, the Hilbert scheme
parametrizing length-$d$ subschemes of $S$, and $\cZ$ is the cycle of the
tautological subscheme of $S\tm S^{[d]}$. One verifies~\cite{Be2,
Prop.~(5)} that there exists $\o^{[d]}\in\G(\O^2_{S^{[d]}})$ which
restricted to the open subset parametrizing reduced subschemes  equals
the push-forward of $p^*\o|_{\cZ}$. Letting $\vf_i\cl B\cdots>S^{[d_i]}$ be
the rational map induced by $\cZ_i$, we can define the terms
appearing in~(2.1) by settting
$$q_{i,*}(p^*\o|_{\cZ_i}):=\vf_i^*\o^{[d_i]}.$$
Mumford~\cite{Mm1} proved that if  $\cZ'\in Z^2(S\tm B)$ is a cycle such
that $\cZ'\cdot S_b\sim \cZ\cdot S_b$ ("$\sim$" denotes rational
equivalence) for all $b\in B$  then $\o_{\cZ'}=\o_{\cZ}$. In particular
we get a well-defined two-form $\o_{\cF}$ on $B$ if we set
$$\o_{\cF}:=\o_{\cZ},\qquad \text{$\cZ\in Z^2(S\tm B)$ a representative
of $c_2(\cF)\in A^2(S\tm B)$.}$$
If $L$ is a line-bundle on $B$ and $\cF':=\cF\ot q^*L$, then  
$\o_{\cF'}=\o_{\cF}$. This allows us to define a two-form $\o_{\xi}$
on the locus $\mx^0\ss\mx$ of smooth (for the reduced structure) stable
points.  More explicitly: if $\mx$ is a fine moduli space we set
$\o_{\xi}:=\o_{\cF}$, where $\cF$ is any tautological family of sheaves
on $S$ parametrized by $\mx$ (the two-form is independent of the choice
of $\cF$), if $\mx$ is not a fine moduli space  one can use a
quasi-tautological family~\cite{Mk2, p.~407} parametrized by $\mx^0$, or
resort to a patching argument. In this section we will deal with the
following question:   Let $[F]\in\mx$ be a generic point, and view
$\o_{\xi}([F])$ as a (skew-symmetric) linear map 
$$\o_{\xi}([F])\cl T_{[F]}\mx \to \O_{[F]}\mx,$$
what is the corank of $\o_{\xi}([F])$? In particular, when does  there exist
an open dense subset of $\mx$ over which $\o_{\xi}$ is a symplectic form?
Before giving a (partial) answer, we must open a digression. 
\definition{(2.2) Definition}
Let $C$ be a smooth irreducible projective curve, and $\t$ be a
theta-characteristic on $C$. We set
$$\l_C(\t,r,d):=h^0(\End_0(V)\ot\t),$$
where $V$ is  the generic stable rank-$r$ vector-bundle on $C$ with $\deg
V=d$. (If $C$ has genus zero then $\l_C$ is not defined, if   $C$ has genus
one $\l_C$ is only defined for $r$, $d$ coprime.)
\enddefinition
A result of Mumford determines the parity of $\l_C$.
\proclaim{(2.3) Proposition~\cite{Mm2}}
Let $\t$ be a theta-characteristic on $C$, and $V$ be a vector-bundle
on $C$. Then
$$h^0\left(\End_0(V)\ot\t\right)\equiv
\left(\rk(V)-1\right)\cdot\left(h^0(\t)+\deg V\right)\pmod{2}.$$
\endproclaim
\demo{Proof}
By Mumford~\cite{Mm2} the quantity $h^0\left(\End_0(V)\ot\t\right)$ is
constant modulo two when $V$ varies in a connected (flat) family.
Since any two vector-bundles on $C$ with the same rank and degree belong
to a connected family, it suffices to check the equation
 for a direct sum of line-bundles; the  computation is left to the
reader. 
\qed
\enddemo
In particular we get
$$\l_C(\t,r,d)\equiv (r-1)\cdot\left(h^0(\t)+d\right)\pmod{2}.$$
\proclaim{(2.4) Conjecture} 
Let $C$ be a smooth irreducible projective curve of genus at least one, and
$\t$ be a theta-characteristic on $C$. Then
$$\l_C(\t,r,d)=\cases
0 & \text{if $(r-1)\cdot\left(h^0(\t)+d\right)\equiv 0 \pmod{2}$,} \\
1 & \text{if $(r-1)\cdot\left(h^0(\t)+d\right)\equiv 1 \pmod{2}$.} 
\endcases$$
\endproclaim
In genus  one the  conjecture  is easily settled, but for bigger
genus we do not know the answer  in general.  When the rank
is two~(2.4) has been proved: in fact there exists a very quick
proof~\cite{L2},  a "Prym variety" proof~\cite{Be1}, and a computational
one~\cite{O1}. Unfortunately we have not succeeded in generalizing any  of
these proofs to higher rank. A different approach, explained at the end of this
section, gives the following.
\proclaim{(2.5) Proposition}
Keep notation as above. Let $C$ be a smooth irreducible projective curve of
genus at least two, and let $\t$ be a theta-characteristic on $C$. Then
$$\l_C(\t,r,h^0(\t))=0.$$
\endproclaim
Now let's go back to moduli of vector-bundles on surfaces.  We will prove
(see~\cite{Mk1,O1}) the following
\proclaim{(2.6) Theorem}
Given  a  polarized surface $(S,H)$ there exists  $\D(r)$ such that the
following holds. Let
$\o$ be a holomorphic two-form on $S$ whose zero-locus $C$ is either
empty or a smooth irreducible curve of genus at least one.  Let
$\xi$ be a set of  sheaf data with $\dx>\D\left(\rk(\xi)\right)$ and, in case
$C$ has genus one, assume also that $\rk(\xi)$, $\left(c_1(\xi)\cdot
K_S\right)$ are coprime. Then the corank of $\o_{\xi}$ at the generic
point of $\mx$ equals  
$$\l_C(K_S|_C,\rk(\xi),c_1(\xi)\cdot K_S).$$
(By convention we set $\l_C=0$ if $C$ is empty.) 
\endproclaim
\remark{Remark}
The lower bound $\D(r)$ of the above theorem is not less than the
quantity $\D(r)$ of Theorem~(0.1), and hence 
$\dim\mx=2r\D_{\xi}-(r^2-1)\chi(\cO_S)$ (here $r:=\rk(\xi)$). A
computation shows that 
$$2r\D_{\xi}-(r^2-1)\chi(\cO_S)
\equiv(r-1)\cdot\left(h^0(K_S|_C)+c_1(\xi)\cdot K_S\right)\pmod{2}.$$
Hence if~(2.4)  holds, Theorem~(2.6) gives that $\o_{\xi}$ is
generically symplectic if $\dim\mx$ is even, and "almost symplectic" if
$\dim\mx$ is odd.
\endremark
Since~(2.4) is true  if the  rank is two, or if $d=h^0(\t)$ by
Proposition~(2.5), we get the following
corollary (see~\cite{Mk1,O1} for the rank-two case) of Theorem~(2.6).
\proclaim{(2.7) Corollary}
Let $(S,H)$ be a  polarized surface and suppose there exists
$\o\in\G(\O^2_S)$ whose zero-locus $C$ is either empty or a smooth
irreducible curve of genus at least one. Let $\xi$ be a set of sheaf data
such that $\dx>\D(\rk(\xi))$, and such that 
$$c_1(\xi)\cdot K_S\equiv h^0(K_S|_C) \pmod{\rk(\xi)}.$$
Then $\o_{\xi}$ is generically a symplectic form. (If $C=\es$ and
$\mx=\mx^{st}$ then $\o_{\xi}$ is a symplectic form~\cite{Mk1}.)
\endproclaim
By the theorem of Mumford on zero-cicles modulo rational
equivalence~\cite{Mm1}  we also get the following.
\proclaim{Corollary}
Let hypotheses be as in the previous corollary. The image of the map 
$$\matrix
\mx & \lra & A^2(S) \\
[F] & \mapsto & c_2(F) 
\endmatrix$$
has dimension equal to $\dim\mx$.
\endproclaim
\remark{(2.8) An example}
Let $\pi\cl S \to \PP^2$ be a double cover
of $\PP^2$ branched over a smooth curve of degree $8n$. A set of
sheaf data $\xi$ with $\det(\xi)=\pi^*\left(\cO_{\PP^2}(n)\right)$ satisfies
the hypotheses of Corollary~(2.7) (for $\dx\gg 0$), hence $\o_{\xi}$
is generically non-degenerate. 
\endremark
\subhead
Proof of Theorem~(2.6)
\endsubhead
We maintain the notation of the introduction to this section. The first step
of the proof consists in identifying $\o_{\cF}$ (up to multiples) with a
certain two-form $\wh{\o}_{\cF}$ introduced by Mukai and
Tyurin~\cite{Mk1,T}. Let $b\in B$, and let 
$$\k\cl T_b(B)\to \Ext^1(\cF_b,\cF_b)$$
be the Kodaira-Spencer map of the family $\cF$. We define
$\wh{\o}_{\cF}$ at $b$ by setting
$$\wh{\o}_{\cF}(v\wedge w):=
\int_S\Tr\left(\k(v)\cup\k(w)\right)\wedge \o.$$
Here "$\cup$" denotes Yoneda pairing, and we are viewing  
$\Tr\left(\k(v)\cup\k(w)\right)$ as a $(0,2)$-form via Dolbeault's
isomorphism. If $\cF_b$ is locally-free then 
$$\Ext^1(\cF_b,\cF_b)\cong H^{0,1}(\End\cF_b),$$
 and $\Tr(\cdot)$ is
obtained   composing  the $(0,1)$-valued endomorphisms $\k(v)$, $\k(w)$,
and taking the trace. A local computation shows that the trace is 
skew-symmetric in this case. For skew-symmetry when $\cF_b$ is not
locally-free see~[M,O1].   
\proclaim{Proposition}
Let notation be as above. Assume the isomorphism class of $\det\cF_b$ is
independent of $b\in B$. Then
$$\o_{\cF}=\left({i\over 2\pi}\right)^2 \wh{\o}_\cF.\tag 2-9$$
\endproclaim
\demo{Sketch of the proof}
First one verifies the following:
\roster
\item
Suppose the isomorphism class of $\det\cF_b$ is
independent of $b\in B$. If $L$ is line-bundle on $S$ and $\cF'=\cF\ot p^*L$,
then 
$$\o_{\cF'}=\o_{\cF}\qquad \wh{\o}_{\cF'}=\wh{\o}_{\cF}.$$
\item
Let 
$$0\to \cE \to \cF \to \cG \to 0$$
be an exact sequence, where $\cE$, $\cF$, $\cG$ are families of torsion-free
sheaves on $S$  with $\det\cE_b$, $\det\cF_b$,
$\det\cG_b$ constant up to isomorphism. Then
$$\o_{\cF}=\o_{\cE}+\o_{\cG} \qquad 
\wh{\o}_{\cF}=\wh{\o}_{\cE}+\wh{\o}_{\cG}.$$
\endroster
Now let's proceed with the proof of~(2.9).  Replacing $B$ by an open
dense subset we can assume there is an exact sequence  
$$0\to\cO_{S\tm B}^{(r-1)}\to \cF\ot p^*\cO_S(nH) \to I_{\cZ}\to 0,$$
where $\cZ$ is a family of zero-dimensional subschemes of $S$
parametrized by $B$, $I_{\cZ}$ is its ideal sheaf,   and $r$ is the rank of
$\cF$. By Items~(1)-(2) above it suffices to prove~(2.9) for
$\cF=I_{\cZ}$. For this it  is enough to consider the universal case:
$B=S^{[d]}$ and $\cZ$ the tautological subscheme of $S\tm S^{[d]}$. There
exists a short locally-free resolution 
$$0\to \cF^1\to\cF^0\to I_{\cZ}\to 0,$$
 where the isomorphism class of $\det\cF^i_{x}$ is independent of $x\in
S^{[d]}$. By Item~(2) it suffices to prove~(2.9) for $\cF^i$. In the de
Rham cohomology of $S^{[d]}$ we have
$$[\o_{\cF^i}]=q_*[c_2(\cF^i)\wedge p_S^*\o],$$
 where $q\cl S\tm S^{[d]}\to S^{[d]}$ is projection (here $c_2(\cF^i)\in
H^4(S\tm S^{[d]})$). On the other hand, by Chern-Weyl theory one
gets~\cite{O1} 
$$\left({i\over 2\pi}\right)^2[\wh{\o}_{\cF^i}]=
q_*[c_2(\cF^i)\wedge p_S^*\o].$$
Hence the two sides of~(2.9) are cohomologous; since they are
holomorphic and since $S^{[d]}$ is projective we conclude that they must be
equal. 
\qed
\enddemo
Now we can prove Theorem~(2.6). By Theorem~(1.1) there exists
$\D(r)$ such that if $\dx>\D(r)$ (where $\rk(\xi)=r$) then
$$\dim W_{\xi}^{2K}  <\text{exp.dim.}\mx. \tag 2-10$$
Furthermore we can assume the generic point on every irreducible
component of $\mx$ represents a stable locally-free sheaf~\cite{O1}. 
Let $[F]\in \mx$ be a generic point; thus $F$ is locally-free, stable, and
by~(2.10) the moduli space is smooth of the expected dimension at  $[F]$
(since $H^0(K)$ has a section, $W^K_{\xi}\ss W^{2K}_{\xi}$). We have
$T_{[F]}\mx\cong H^1(\End_0(F))$ (see Secion~(1)), and by~(2.9) 
$$\o_{\xi}(v\wedge w)=
\left({i\over 2\pi}\right)^2\int_S\Tr\left(v\wedge(w\cdot\o)\right).$$
By Serre duality the bilinear map
$$\matrix
H^1(\End_0(F))\tm H^1(\End_0(F)\ot K) & \lra & H^2(K)\cong\CC \\
(\a,\b) & \mapsto &\Tr\left(\a\cup\b\right)
\endmatrix$$
is a perfect pairing, hence it suffices to show that for generic $[F]\in \mx$ 
the map
$$H^1(\End_0(F))\brel \ot\o\over\lra H^1(\End_0(F)\ot K)$$
has corank $\l_C(K_S|_C,\rk(\xi),c_1(\xi)\cdot K_S)$. This certainly holds if
$C=\es$, hence we can assume $C\not=\es$. Consider the exact sequence  
$$H^0(\End_0(F)\ot K)\to H^0(\End_0(F)\ot K|_C)\to 
H^1(\End_0(F))\brel \ot\o\over\lra H^1(\End_0(F)\ot K).$$
Since $[F]$ is generic and since~(2.10) holds, we have $h^0(\End_0(F)\ot
K)=0$. Thus we must show that
$$h^0(\End_0(F)\ot K|_C)=
\l_C(K_S|_C,\rk(\xi),c_1(\xi)\cdot K_S).\tag 2-11$$
Hence it suffices to prove that if $[F]\in\mx$ is generic
then $F|_C$ is the generic stable vector bundle (of rank $\rk(\xi)$ and
determinant $\det(\xi)|_C$). So let $[E]\in \mx$ with $E$ locally-free, stable
and $[E]\notin W^{2K}_{\xi}$; we claim the map $\rho\cl
\Def^0(F)\to\Def^0(F|_C)$ defined by restriction is surjective. In fact
both the domain and  codomain are smooth, and the differential  $d\rho$
fits into the exact sequence 
$$\multline
H^1(\End_0(E))\brel d\rho\over\lra H^1(\End_0(E)|_C)\to\\ 
\to H^2(\End_0(E)\ot[-K])\cong H^0(\End_0(E)\ot[2K])^*=0.
\endmultline$$ 
By hypothesis $C$ has genus at least two or, if it has genus  one,  $\rk(\xi)$
and $c_1\left(\det(\xi)|_C\right)$ are coprime. Thus the generic
vector-bundle parametrized by $\Def^0(E|_C)$ is the generic stable bundle
with the chosen rank and determinant. By surjectivity of $\rho$ we
conclude that Equation~(2.11) holds for the generic $[F]\in X$.
\subhead
Proof of~(2.5)
\endsubhead
We let $C$ be a smooth irreducible projective curve of genus at least
two. We will examine vector-bundles $\cE$ obtained as extensions 
$$0\to V^{*}\to\cE\to\cO_C\to 0.\tag 2-12$$
\proclaim{(2.13) Lemma}
Keep notation as above. Assume that:
\roster
\item
$h^0(\End_0(V)\ot\t)=0$.
\item
$0\le\deg V\le h^0(\t)$.
\item
$V$ is generic among vector-bundles of the same degree and rank.
\item
The extension class $\eta\in H^1(V^{*})$ of~(2.12) is generic.
\endroster
Then there is a natural identification
$$H^0\left(\End_0(\cE)\ot\t\right)\cong
\Ker\left(H^0(\t)\brel\del\over\lra H^1(V^*\ot\t)\right),$$
where $\del$ is the coboundary map associated to the sequence
obtained  tensoring~(2.12) by $\t$:
$$0\to V^{*}\ot\t\to\cE\ot\t\to\cO_C(\t)\to 0.\tag 2.14$$
\endproclaim
\demo{Proof}
Scalar endomorphisms give an injection $H^0(\t)\hra
H^0\left(\End(\cE)\ot\t\right)$, and there is a splitting 
$$H^0\left(\End(\cE)\ot\t\right)= H^0(\t)\op
H^0\left(\End_0(\cE)\ot\t\right).$$
Thus it suffices to give an identification
$$H^0\left(\End(\cE)\ot\t\right)/H^0(\t)\cong
\Ker\left(H^0(\t)\brel\del\over\lra H^1(V^*\ot\t)\right).$$
Let $\vf\in H^0\left(\End(\cE)\ot\t\right)$. First we prove 
$$\vf(V^{*})\ss V^{*}\ot\t.\tag{$*$}$$
For this it suffices to show that
$\Hom(V^{*},V^{*}\ot\t)\hra \Hom(V^{*},\cE\ot\t)$ is an isomorphism.
Tensoring~(2.14) by $V$ we get a coboundary map
$$\matrix
H^0(V\ot\t) &\brel\del_V\over\lra & H^1(V\ot V^{*}\ot\t). \\
\a & \mapsto & \a \cup \eta 
\endmatrix$$
We must show $\del_V$ is injective. Consider the trace map
$$H^1(V^{*}\ot V\ot\t)\brel\Tr\over\to H^1(\t).$$
We will prove that $\Tr\circ\del_V$ is injective. By Serre duality we can
view $\Tr\circ\del_V$ as a map 
$$\Tr\circ\del_V\cl H^0(V\ot\t) \to H^0(\t)^{*}.$$
Explicitly, since the extension class  $\eta$ of~(2.12) is an element of
$H^0(V\ot K_C)^{*}$ (by Serre duality), we have 
$$\langle\Tr\circ\del_V(\a),\b\rangle=\langle\eta,\a\ot\b\rangle,
\qquad \a\in H^0(V\ot\t),\ \b\in H^0(\t).$$
Because the map 
$$\matrix
H^0(\t) & \lra & H^0(V\ot K_C) \\
\b & \mapsto & \a\ot\b
\endmatrix$$
is injective for any non-zero $\a\in H^0(V\ot\t)$,  there is a
well-defined map 
$$\matrix
\PP:=\PP\left(H^0(V\ot\t)\right) & \brel\Phi\over\lra & 
Gr:=Gr\left(h^0(\t),h^0(V\ot K_C)\right), \\
[\a] & \mapsto & \{\a\ot\b|\ \b\in H^0(\t)\}
\endmatrix$$
where $Gr(m,n)$ is the Grassmannian of $m$-dimensional vector subspaces
of $\CC^n$. Let
$$\L:=\bigcup_{[\a]\in\PP} 
\{\text{$[\eta]\in \PP\left(H^0(V\ot K_C)^{*}\right)$  vanishes on
$\Phi([\a])$}\}.$$     
(Notice that $H^0(V\ot K_C)\not=0$ because $\deg V\ge 0$ by Item~(2), and
because $C$ has genus at least two.) We must show that
$$\L\not=\PP\left(H^0(V\ot K_C)^{*}\right).\tag{$\bu$}$$
A dimension count gives
$$\aligned
\dim\L & \le h^0(V\ot\t)-1+h^0(V\ot K_C)-h^0(\t)-1\\
& =\dim\PP\left(H^0(V\ot K_C)^{*}\right)
-\left(h^0(\t)-h^0(V\ot\t)+1\right). 
\endaligned\tag{$\dag$}$$
By our hypotheses  $\deg V\ge 0$ and $V$ is generic. This  implies that
$$h^0(V^{*}\ot\t)=0,\tag{$2.15$}$$
and thus by Serre duality $h^1(V\ot\t)=0$.
Hence 
$$h^0(V\ot\t)=\chi(V\ot\t)=\deg V.\tag{$2.16$}$$
By~($\dag$) we conclude that if  $\deg V\le h^0(\t)$ then~($\bu$) holds.
Thus for $\eta$  generic the map $\del_V$ is injective. (Of course $\deg V\le
h^0(\t)$ is also necessary for $\del_V$ to be injective.) This proves~($*$).
Now we can finish the proof of the lemma. By~($*$) and Item~(1) the
restriction of $\vf$ to $V^{*}$ is equal to scalar multiplication by a certain
section $\s\in H^0(\t)$; thus 
$$\left(\vf-\s\right)(V^{*})=0,$$
or in other words
$$H^0(\End(\cE)\ot\t)/H^0(\t)\cong
\Hom\left(\cO_C,\cE\ot\t\right)=H^0(\cE\ot\t).$$ 
Writing out the long exact cohomolgy
sequence associated to~(2.14),  the lemma  follows from~(2.15).
\qed
\enddemo
Let's prove Proposition~(2.5) by induction on the
rank $r$. The case $r=1$ is trivial. Let's prove the inductive step.
We assume that $V$ is the generic stable rank-$r$ vector-bundle with $\deg
V=h^0(\t)$, and that~(2.5) holds for $V$. Consider the generic
extension~(2.12). If we show that $h^0(\End_0(\cE)\ot\t)=0$ then we are
done, because $\rk(\cE^{*})=(r+1)$ and $\deg\cE^{*}=h^0(\t)$. By 
Lemma~(2.13) it suffices to prove that 
$$\matrix
H^0(\t) & \brel\del\over\lra & H^1(V^{*}\ot\t) \\
\b & \mapsto & \b\cup\eta 
\endmatrix$$
is injective. This coboundary is the transpose of the map $\Tr\circ\del_V$
appearing in the proof of Lemma~(2.13).  We have proved
$\Tr\circ\del_V$ is injective, and thus $\del$ is injective if and only if
$h^0(\t)=h^1(V^{*}\ot\t)$. By Serre duality $h^1(V^{*}\ot\t)=h^0(V\ot\t)$,
hence the result follows from~(2.16).
\head
3.  Kodaira dimension of the moduli space.
\endhead
The main result is due to Jun Li~\cite{L2}: moduli spaces of
rank-two vector-bundles on a surface of general type are often of general
type. We will prove Jun Li's theorem for fine moduli spaces, with some
additional hypotheses. The proof in general is more difficult, the main
problem being the  analysis of  singularities coming from strictly
semistable sheaves~\cite{L2}.  
 Jun Li's theorem, for fine moduli spaces,  extends to higher rank if
Conjecture~(2.4) is true. In particular by Proposition~(2.5) many
higher-rank moduli spaces are of general type; we will give some examples.
We will also quickly mention some results concerning moduli spaces on
surfaces not of general type. In the proofs we will usually choose a
particularly nice polarization (essentially a multiple of $K$): this is not a
significant restriction because of Remark~(0.2).

Throughout this section we assume $\mx$ is a fine moduli space; we let
$\cF$ be a tautological sheaf on $S\tm\mx$. To simplify notation we set
$r:=\rk(\xi)$.  
\subhead
The canonical line-bundle
\endsubhead
Let $\mx^{sm}$ be the subscheme of $\mx$ parametrizing stable sheaves
$F$ with vanishing obstruction space, i.e.~such that
$$\Ext^2(F,F)^0=0.$$
By deformation theory $\mx^{sm}$ is smooth. 
\proclaim{(3.1) Lemma}
Keep notation as above. Then modulo torsion 
$$K_{\mx^{sm}}\cong\cL(rK_S).$$
\endproclaim
\demo{Proof}
Let $\pi\cl S\tm\mx^{sm}\to\mx^{sm}$ be the projection, and let $\cF^{sm}$
be the restriction of $\cF$ to $S\tm\mx^{sm}$. Define
$Ext^p_{\pi}(\cF^{sm},\cF^{sm})^0$ as the sheaf on $\mx^{sm}$ fitting into
the exact sequence 
$$0 \to Ext^p_{\pi}(\cF^{sm},\cF^{sm})^0 \to  Ext^p_{\pi}(\cF^{sm},\cF^{sm}) 
\brel Tr \over \lra R^p\pi_*\cO \to 0.\tag 3-2$$
Since $Ext^p_{\pi}(\cF^{sm},\cF^{sm})$ is a vector-bundle with fiber 
$\Ext^p(F,F)$ over $[F]\in\mx^{sm}$, the fiber of
$Ext^p_{\pi}(\cF^{sm},\cF^{sm})^0$ over  $[F]$ is canonically isomorphic to
$Ext^p(F,F)^0$. Thus by deformation theory 
$$T\mx^{sm}\cong Ext^1_{\pi}(\cF^{sm},\cF^{sm})^0.$$
Exact sequence~(3.2) for $p=1$ gives
$$c_1\left(Ext^1_{\pi}(\cF^{sm},\cF^{sm})^0\right)=
c_1\left(Ext^1_{\pi}(\cF^{sm},\cF^{sm})\right).$$
On the other hand, by definition of $\mx^{sm}$ we have 
$Ext^p_{\pi}(\cF^{sm},\cF^{sm})^0=0$ for $p=0,2$, and hence~(3.2) 
gives
$$c_1\left(Ext^p_{\pi}(\cF^{sm},\cF^{sm})\right)=0\qquad p=0,2.$$
>From the above equalities we get that  in the Chow group
$A^1(\mx^{sm})_{\QQ}$ 
$$c_1(K_{\mx^{sm}})=
-c_1\left(Ext^1_{\pi}(\cF^{sm},\cF^{sm})^0\right)=\sum_{p=0}^2(-1)^p
c_1\left(Ext^p_{\pi}(\cF^{sm},\cF^{sm})\right).\tag 3-3$$ 
The right-hand side can be computed by applying
Grothendieck-Riemann-Roch: setting 
$$ch(\cF^{sm})^*:=\sum_n (-1)^nch_n(\cF^{sm}),$$ 
we have
$$\align
\sum_{p=0}^2(-1)^p c_1\left(Ext^p_{\pi}(\cF^{sm},\cF^{sm})\right)= &
\pi_{*}\left[ch(\cF^{sm})^*ch(\cF^{sm})\text{Td}(S)\right]_3 \\ 
=&\pi_{*}\left[\left(c_2(\cF^{sm})- {r-1\over
2r}c_1(\cF^{sm})^2\right)\cdot rK_S\right]. 
\endalign$$  
The lemma  follows from the above formula together with~(3.3)
and~(1.12). 
\qed
\enddemo
\subhead
Surfaces of Kodaira dimension at most one
\endsubhead
First assume $S$ is a Del Pezzo surface, and let $H=-K$. Since $K\cdot H<0$
the obstruction space
$$\Ext^2(F,F)^0\cong\left(\Hom(F,F\ot K)^0\right)^{*}$$
vanishes for every $[F]\in\mx$, hence $\mx$ is smooth (of the expected
dimension). Thus Lemma~(3.1), together with Theorem~(1.13), implies
that $\k(\mx)=-\infty$. In fact more is known~\cite{ES}: if $S=\PP^2$ the
moduli space is often rational. More in general, it is natural to expect that
if $S$ is (birationally) ruled, the moduli space is uniruled (for $\dx>>0$).
Hoppe and Spindler~\cite{HS} treat the case of rank two. 

If $S$ is a $K3$ surface, the moduli space is smooth, hence by Lemma~(3.1)
we conclude that $\k(\mx)=0$. In fact more is true: if $\o$ is a
non-zero two-form on $S$, the  two-form $\o_{\xi}$ is everywhere
non-degenerate~\cite{Mk1}, thus $\mx$ is holomorphically symplectic. 

Finally let's consider the case of a minimal surface of Kodaira dimension
one. Let
$$f\cl S\to B$$
be the elliptic fibration, and for $b\in B$ let $C_b:=f^*(b)$. We assume
that the set of sheaf data satisfies:
$$\text{$\rk(\xi)$ and $c_1(\xi)\cdot C_b$ are coprime.}$$
It is convenient to choose the polarization $H$ to be very close to $C_b$
in the N\'eron-Severi group $NS(S)$ (how close will depend on $\rk(\xi)$
and $\dx$). Such a polarization is called {\it suitable}~\cite{F1}. 
\proclaim{(3.4) Lemma}
Let notation and hypotheses be as above. Then there are no strictly
$H$-slope-semistable torsion-free sheaves on $S$. Furthermore, a
torsion-free sheaf $F$ on $S$ is $H$-slope-stable if and only if $F|_{C_b}$
is stable for the generic elliptic fiber $C_b$. 
\endproclaim
Thus $\mx^{st}=\mx$; furthermore $\mx$ is a fine moduli space (apply
Remark~(A.7) of~\cite{Mk2}). One also verifies that $\mx$ is
smooth~\cite{F3}, hence $\k(\mx)$ equals the dimension of
$$X:=\text{\bf Proj}
\left(\bigoplus_{n=0}^{\infty}H^0(K_{\mx}^{\ot n})\right).$$ 
We expect that $X$ and the canonical map $\mx\to X$ are described as
follows, but we have not checked the details. A computation shows that
$\dim\mx$ is even, so set $\dim\mx=2n$. If $[F]\in\mx$, then by
Lemma~(3.4) the restriction to the generic elliptic fiber $C_b$ is stable,
but there are $n$ fibers $C_{b_1},\ldots,C_{b_n}$ such that $F|_{C_{b_i}}$
is not stable (or not locally-free). Thus we get a morphism
$$\matrix
\mx & \brel \Phi\over \lra & B^{(n)} \\
[F] & \mapsto & b_1+\cdots + b_n. 
\endmatrix$$
Then the canonical model $X$ is identified with $B^{(n)}$, and the canonical
map is identified with $\Phi$. Thus
$${\k(\mx)\over\dim(\mx)}={1\over 2}={\k(S)\over\dim(S)}.$$
\subhead
Surfaces of general type
\endsubhead
We will prove the following result.
\proclaim{(3.5) Theorem (Jun Li~\cite{L2})}
Let $S$ be a  surface with ample canonical bundle,  and let $H$ be a rational
multiple of $K_S$. Let $\xi$ be a set of sheaf data on $S$ such that: 
\roster
\item
The moduli space $\mx$ is fine.
\item
The codimension of $W_{\xi}^K$ in $\mx$ is at least two.
\item
There exists $\o\in\G(\O^2_S)$ such that $\o_{\xi}$ is generically
non-degenerate.
\endroster
Then $\mx$ is of general type.
\endproclaim
Combining the above theorem with~(1.1), (2.6) and~(2.7) one
gets the following corollaries.
\proclaim{Corollary}
Let $S$, $H$ be as in the statement of Theorem~(3.5), and assume there
exists $\o\in\G(\O^2_S)$ whose zero-locus is a smooth irreducible curve
(canonical) curve $C$. There exists a function $\D(r)$ for which the
following holds. Let $\xi$ be a set of sheaf data such that:
\roster
\item
$\mx$ is a fine moduli space .
\item
$\dx>\D(\rk(\xi))$.
\item
$\l_C\left(K_S|_{C},\rk(\xi),c_1(\xi)\cdot C\right)=0$.
\endroster
Then $\mx$ is of general type.
\endproclaim
\proclaim{Corollary}
Let hypotheses be as above, except that we replace Item~(3) by: 
\roster
\item[4]
$c_1(\xi)\cdot C\equiv h^0(K_S|_C) \pmod{\rk(\xi)}$. 
\endroster
Then $\mx$ is of general type.
\endproclaim
\remark{(3.6) An example}
Let $\pi\cl S \to \PP^2$ be a double cover branched over a smooth
curve of degree $8n$, and let $H:=\pi^*\cO_{\PP^2}(1)$. Let $\xi$ be a set of
sheaf data such that $\det(\xi)=nH$, $\dx>\D(\rk(\xi))$, and 
$$\gcd\left\{\rk(\xi),2n,n^2+n-c_2(\xi)\right\}=1.$$ 
(This last condition ensures that $\mx$ is a fine moduli
space~\cite{Ma,Mk2}.)   Then the hypotheses of the corollary are satisfied,
hence $\mx$ is of general type. 
\endremark
Let's prove Theorem~(3.5). By Items~(1)-(2), together with deformation
theory, the moduli space is a local complete intersection, hence its
dualizing sheaf is a line-bundle; we denote it by $K_{\mx}$. By Item~(2)
$\mx$ is smooth in codimension one, hence Lemma~(3.1) gives that
$$K_{\mx}\cong\cL(rK) \qquad \text{up to torsion.}$$
Since $K$ is a positive multiple of $H$ the (fractional) line-bundle
$\cL(rK)$ is big, by Theorem~(1.13).  Hence there exists a positive $c$ such
that for $n$ large enough and  sufficiently divisible   
$$\G(K_{\mx}^{\ot n})=cn^d+O(n^{d-1}), \tag 3-7$$
where $d:=\dim\mx$. This is not sufficient to
conclude that $\mx$ is of general type, because of the
presence of singularities.  Let $\rho\cl\wt{\cM}_{\xi}\to\mx$ be a
desingularization. We have   
$$\rho^*K_{\mx}=K_{\wt{\cM}_{\xi}}(\sum_i a_i E_i),$$
for some $a_i\in\ZZ$, where $E_i$ are the exceptional divisors of $\rho$.
Let $a$ be a non-negative number such that $a\ge a_i$ for all $i$. By the
above equation we have
$$\rho^*\G\left(K_{\mx}^{\ot n}\right)\ss
\G\left(K_{\wt{\cM}_{\xi}}^{\ot n}(anE)\right),$$
where $E:=\sum_i E_i$. Now we will use the two-form $\o_{\xi}$ to
produce  a non-zero section of $K_{\wt{\cM}_{\xi}}(-E)$. Let
$\wt{\cF}$ be the bull-back to $S\tm\wt{\cM}_{\xi}$ of the tautological
family over $S\tm\mx$; then 
$$\s:=\wedge^{d/2}\o_{\wt{\cF}}\in\G\left(K_{\wt{\cM}_{\xi}}\right)$$
is non-zero because $\o_{\xi}$ is generically non-degenerate, by
Item~(3). Since $\dim\rho(E_i)<\dim E_i$, the Kodaira-Spencer map of
$\wt{\cF}$ has a non-trivial kernel at the generic point of $E_i$, hence by
Equation~(2.9) the  two-form $\o_{\wt{\cF}}$ is degenerate along
$E_i$; thus 
$$\s\in \G\left(K_{\wt{\cM}_{\xi}}(-E)\right)$$
At this point we are done: letting $N:=n(a+1)$, 
we have an injection 
$$\matrix
\rho^*\G\left(K_{\mx}^{\ot n}\right) & \hra &
\G\left(K_{\wt{\cM}_{\xi}}^{\ot N}\right) \\
\tau & \mapsto & \tau \ot \s^{\ot(an)}. 
\endmatrix$$
By~(3.7) we conclude that $\wt{\cM}_{\xi}$ is of general type, hence
$\mx$ is of general type. 

\enddocument